\documentclass[11pt]{amsart}

\newtheorem{theorem}{Theorem}[section]

\newtheorem{remark}[theorem]{Remark}
\newtheorem{proposition}[theorem]{Proposition}
\newtheorem{corollary}[theorem]{Corollary}
\newtheorem{problem}[theorem]{Problem}
\theoremstyle{definition}

\newtheorem{example}[theorem]{Example}
\numberwithin{equation}{section}

\def\a{{\bf a}}
\def\b{{\bf b}}

\def\d{{\bf d}}

\def\x{{\bf x}}
\def\y{{\bf y}}
\def\1{{\bf 1}}
\def\0{{\bf 0}}
\def\IC{{\mathbb C}}
\def\IR{{\mathbb R}}

\def\diag{{\rm diag}\,}

\def\tr{{\rm tr}\,}

\def\diag{{\rm diag}\,}
\def\[{\left [}
\def\]{\right ]}
\def\({\left (}
\def\){\right )}
\def\la{{\langle}}
\def\ra{{\rangle}}
\def\Ra{{\ \Rightarrow\ }}

\def\Lra{{\ \Leftrightarrow\ }}
\def\dfrac{\displaystyle\frac}

\begin{document}
\openup 1.1 \jot

\title[Interpolation Problems by Completely positive maps]
{Interpolation Problems by Completely positive maps}

\author{Chi-Kwong Li}
\address{Department of Mathematics, College of William \& Mary,
Williamsburg, VA 23187}
\email{ckli@math.wm.edu}

\author{Yiu-Tung Poon}
\address{Department of Mathematics, Iowa State University,
Ames, IA 50051}
\email{ytpoon@iastate.edu}

%\author{Nung-Sing Sze}
%\address{Department of Mathematics, University of Connecticut,
%Storrs, CT 06269}
%\email{sze@math.uconn.edu}
%\thanks{}


%\date{}

\dedicatory{}

%    "Communicated by" -- provide editor's name; required.
\commby{}

\begin{abstract}
Given commuting families  of Hermitian matrices $\{A_1, \dots, A_k\}$ and
$\{B_1, \dots, B_k\}$,
conditions for the existence of a completely positive map $\Phi$,
such that $\Phi(A_j) = B_j$ for $j = 1, \dots, k$, are studied.
Additional properties such as unital or / and trace preserving
on the map $\Phi$ are also considered.
Connections of the study to dilation theory, matrix inequalities,
unitary orbits, and quantum information science are mentioned.
\end{abstract}

\maketitle

\section{Introduction}

Denote by $M_{n,m}$ the set of $n\times m$ complex matrices,
and use $M_n$ to denote $M_{n,n}$.
Let $H_n$  be the set of Hermitian matrices in $M_n$.
A matrix $A \in H_n$ is positive semidefinite if all
eigenvalues of $A$ are nonnegative.
A linear map $\Phi: M_n \rightarrow M_m$ is {\it positive} if it
maps positive semidefinite matrices to positive semidefinite matrices.
Suppose $M_k(M_n)$ is the algebra of block matrices of the form
$(A_{ij})_{1 \le i, j \le k}$ with $A_{ij} \in M_n$ for each
pair of $(i,j)$.
A linear map $\Phi: M_n \rightarrow M_m$ is
{\it completely positive}
if for each positive integer $k$, the map
$I_k \otimes \Phi: M_k(M_n) \rightarrow M_k(M_m)$
defined by $(I_k \otimes \Phi)(A_{ij}) = (\Phi(A_{ij}))$
is positive.

The purpose of this paper is to study the following.

\begin{problem}\label{1.1}
Given  $A_1, \dots, A_k  \in M_n$ and $B_1, \dots, B_k \in M_m$,
determine the necessary and sufficient condition
for the existence of a completely positive map
$\Phi: M_n \rightarrow M_m$ such that $\Phi(A_j) = B_j$
for $j = 1, \dots, k$, and possibly with the
additional properties that
$\Phi(I_n) = I_m$ or/and $\Phi$ is trace preserving.
\end{problem}

\medskip
Clearly, this can be viewed as an interpolation problem
by completely positive maps.
Denote by $H_n$ the set of Hermitian matrices in $M_n$.
Since for a positive linear map $\Phi$ satisfying
$\Phi(X) = Y$ if and only if
$\Phi( X^*) =   Y^*$ for any $(X,Y) \in M_n\times M_m$,
we can focus on the study of Problem \ref{1.1} for
$\{A_1, \dots, A_k\} \subseteq H_n$ and $\{B_1, \dots, B_k\} \subseteq H_m$

\medskip
Over half a century ago, Steinspring \cite{S} introduced completely positive
maps in the study of dilation problems for operators.
Since then, the area has been studied extensively \cite{Pau}.
In particular, researchers have obtained
interesting structure theorem for completely positive maps on
matrices. For example, Choi \cite{C} (see also \cite{K}) showed that
a linear map $\Phi: M_n\to M_m$ is completely positive if and only if there exist
$F_1, \dots, F_r\in M_{n,m}$ such that
\begin{equation}\label{operatorsum}
\Phi(A)=\sum_{j=1}^rF_j^*AF_j.
\end{equation}
This is called an {\it operator sum representation}
of the completely positive map $\Phi$.

\medskip
In the context of quantum information theory,
every quantum operation is a completely positive map
sending quantum states to quantum states, where quantum states are
represented as {\it density matrices}, i.e., positive semidefinite
matrices with trace one. Because a completely positive map
$\Phi: M_n \rightarrow M_m$ representing a quantum operation
will send density matrices to density matrices, the map
$\Phi$ is {\it trace preserving}, i.e., $\Phi(A)\in M_m$ and $A \in M_n$
always have the same trace. Therefore, in quantum information science,
most studies are on trace preserving completely positive maps. On the
other hand, in the $C^*$-algebra context, since the trace function may
not be defined, most studies are on
{\it unital}  maps,  i.e., $\Phi(I_n) = I_m$.
While rich theory has been developed for completely positive
maps, for example, see \cite{Pau}, there are many basic
problems motivated by applied topics which deserve further study.

In quantum information
science, one has to study and construct
quantum operations sending a specific family of
density operators to another family.
This clearly reduces to Problem \ref{1.1}
if we restrict our attention to
density matrices $A_1, \dots, A_k \in M_n$ and
$B_1, \dots, B_k \in M_m$; see \cite{NC}.

\medskip
Using the operator sum representation (\ref{operatorsum})
of completely positive maps, and the
inner product $\langle A, B \rangle = \tr(AB^*)$
for $A, B \in M_{m,n}$, one has
the following result showing the connection between
trace preserving completely positive maps and
unital completely positive maps.

\begin{proposition}\label{1.2}
Suppose $\Phi: M_n \rightarrow M_m$
is a completely positive map
with operator sum representation  in (\ref{operatorsum}).
Then $\Phi$ is unital if and only if
$\sum_{j=1}^rF_j^*F_j=I_m$;
$\Phi$ is trace preserving if and only if
$\sum_{i=1}^rF_jF_j^*=I_n$.
Moreover, the dual linear map $\Phi^*: M_m \rightarrow M_n$
defined by
$$\Phi^*(B) = \sum_{j=1}^r F_jBF_j^*$$
is the unique linear map satisfying
$\la\Phi(A),B\ra = \la A, \Phi^*(B)\ra$ for all  $(A,B) \in M_n \times M_m$.
Consequently, $\Phi$ is unital if and only if $\Phi^*$ is trace preserving.
\end{proposition}

The above proposition provides a link between problems and results for
trace preserving completely positive maps and unital completely positive maps.
Therefore, one might expect that the results and proofs for the two types of problems
can be converted to each other easily. However, this does
not seem to be the case as shown in our results.
In fact, some results in trace preserving completely positive maps
are more involved in our study, and they
have no analogs for unital completely positive maps; see
Remark \ref{rmk2.2} b), and the remarks after Corollary \ref{new3.4} and Theorem \ref{3.5}.

It is known that the study of completely positive maps
are closely related to the dilations of operators.
Recall that a matrix $B\in M_m$ has a {\it dilation} $A\in M_n$ if there is
an $n\times m$ matrix $V$ such that $V^*V  = I_m$ and $V^*AV = B$.
The next result shows that Problem \ref{1.1} can be formulated as
problems involving dilations and principal submatrices of a matrix.

\begin{proposition} \label{1.3}
Suppose $\Phi: M_n \rightarrow M_m$
is a completely positive map
with operator sum representation $(\ref{operatorsum})$.
If $F   = \[\begin{matrix}F_1  \cr \vdots \cr F_r \cr\end{matrix}\]$, then
$\Phi(A) = F^*(I_r \otimes A)F$.
If $\tilde F  = \[ F_1\, \cdots \ F_r\]$ and
$\tilde F^* A \tilde F = (A_{ij})$ with $A_{11}, \dots, A_{rr} \in M_m$, then
$\Phi(A) = A_{11} + \cdots + A_{rr}$.  Furthermore, the following hold.
\begin{itemize}
\item[{\rm (a)}]
The map $\Phi$ is unital if and only if $F^*F = I_m$.
\item[{\rm (b)}]
The map $\Phi$ is trace preserving if and only if
 $\tilde F \tilde F^*= I_n$.
\end{itemize}
\end{proposition}

\begin{proposition} \label{1.4}
Suppose $\Phi: M_n \rightarrow M_m$ is a completely positive map. Given unitaries  $U\in M_n$ and $V\in M_m$, define $\Psi: M_n \rightarrow M_m$ by $\Psi(X)=V^*\Phi(U^*XU)V $. Then $\Psi$ is also completely positive. Furthermore, $\Psi$ is unital and/or trace preserving if and only if $\Phi$ has the corresponding property. If
$\Phi(X)= \sum_{j=1}^r  F_j ^*X  F_j$, then $ \Psi(X)= \sum_{j=1}^r  (UF_jV)^*X(UF_jV)$.
\end{proposition}

Propositions \ref{1.2} -- \ref{1.4} will be used in the
subsequent discussion. Our paper is organized as follows.
In Section 2, we determine the condition
for the existence of
completely positive maps (possibly with additional
conditions such as unital or/and trace presering)
$\Phi: M_n \rightarrow M_m$
sending a given commuting family of matrices
in $H_n$ to another one in $H_m$.
In Section 3, we give a more detailed analysis
for the case when each family has only one matrix.
Some related results, additional remarks, and open problems
will be mentioned in Section 4.

\section{completely positive maps between commuting families}

In this section, we consider (unital) completely positive maps
$\Phi: M_n \rightarrow M_m$ sending a given commuting family
of matrices in $H_n$ to another commuting family in $H_m$.

\begin{theorem} \label{new2.1} Let $\{A_1, \dots, A_k\}\subseteq H_n$ and
$\{B_1, \dots, B_k\} \subseteq H_m$ be two commuting families.
Then there exist unitary matrices  $U \in M_n$ and $V \in M_m$   such that
$U^*A_iU  $ and $VB_iV^* $ are diagonal matrices with diagonals
$\a_i=( a_{i1}, \dots, a_{in})$ and
$\b_i =  (b_{i1}, \dots, b_{im})$ respectively,
for $i = 1, \dots, k$. The following conditions are equivalent.
\begin{itemize}
\item [{\rm (a)}] There is a  completely positive map
$\Phi: M_n \rightarrow M_m$ such that $\Phi(A_i) = B_i$ for
$i=1, \dots, k$.
\item[{\rm (b)}] There is an $n\times m$ nonnegative
matrix $D = (d_{pq})$ such that
$$(b_{ij}) = (a_{ij}) D.$$
\end{itemize}
Suppose {\rm (b)} holds. For $1\le j\le m$, let   $F_j$ be the $n\times n$
matrix  having the $j${\rm th} column equal to $(\sqrt{d_{1j}},\sqrt{d_{2j}},
\dots,\sqrt{d_{nj}})^t$ and zero elsewhere. Then we have
\begin{equation}\label{ba}B_i = \sum_{j=1}^r (UF_jV)^*A_i(UF_jV),
\qquad i = 1, \dots, k,\end{equation}
Furthermore,
\begin{enumerate}
\item $\Phi$ in {\rm (a)} is unital if and only if $D$ in {\rm (b)} can be chosen to be column stochastic.
\item $\Phi$ in {\rm (a)} is trace preserving if and only if $D$ in {\rm (b)}
 can be chosen to be row stochastic.
\item $\Phi$ in {\rm (a)} is unital and trace preserving if and only if $D$ in {\rm (b)}
can be chosen to be doubly stochastic.
\end{enumerate}
\end{theorem}

\it Proof. \rm By Proposition \ref{1.4}, we may assume that $A_i=\diag(\a_i)$ and $B_i=\diag(\b_i)$ for $i=1,\dots,k$, and take $U=I_n$, $V=I_m$ in (\ref{ba}).

(a) $\Ra $ (b):
 Suppose there is a completely positive map $\Phi: M_n \rightarrow M_m$ such that $\Phi(A_i) = B_i$
 for $i=1, \dots, k$.   By (\ref{operatorsum}), we have $F^j=\(f^j_{pq}\)\in M_{n,m}$, $j=1,\dots,r$,
 such that
 $$\Phi(X)=\sum_{j=1}^rF_j^*XF_j.$$
For $1\le p\le n$ and $1\le q\le m$, let $d_{pq}= \sum_{j=1}^r|f^j_{pq}|^2 $. Then $D=(d_{pq})$ is  an $n\times m$ nonnegative
matrix   such that
$$(b_{ij}) = (a_{ij}) D.$$

(b) $\Ra$ (a): Suppose  $D=(d_{pq})$ is a nonnegative matrix satisfying
$(b_{ij})=(a_{ij})D$. Let $F_j$ be defined as in the theorem.
Then direct computation shows that for every $\x=(x_1,\dots,x_n)$,
$\sum_{j=1}^rF_j^*\diag(\x)F_j$ is a diagonal matrix with diagonal $\y=\x D$.
Hence, for $i=1,\dots,k$, we have
$$\diag(\b_i)=\sum_{j=1}^mF_j^*\diag(\a_i)F_j, $$
for $i=1,\dots,k$.

For $N\ge 1$, let $\1_N$ be a row vector of $N$ 1's, then  we have
 $$ \Phi \mbox{  is unital }\Lra \Phi(I_n)=I_m\Lra \1_m=\1_n D
\Lra  D \mbox{  is column stochastic.}$$

This proves (1). The proof for cases (2) and (3) are similar.
\qed

A completely positive map $\Phi:M_n\to M_n$ is called {\it mixed unitary}
if there exist unitary matrices $U_1,\dots,U_r \in M_n $ and
positive numbers $t_1, \dots, t_r$ summing up to 1 such that
$$\Phi(X) = \sum_{j=1}^r  t_j U_j^*X U_j.$$
Clearly, every mixed unitary completely positive map is unital and trace
preserving. For $n\ge 3$, there exists \cite{LS} a unital trace preserving
completely positive map which is not mixed unitary.

By the Birkhoff Theorem \cite{MO}, a doubly stochastic matrix  $D$ can be expressed in the form \begin{equation}\label{ds}D = \sum_{j=1}^r t_j P_j\end{equation}
for some positive numbers $t_1, \dots, t_r$ summing up to 1 and
permutation matrices $P_1, \dots, P_r \in M_n$.
Using this result and (3) in Theorem \ref{new2.1}, we have the following corollary.

\begin{corollary} \label{new2.2}
Under the hypothesis of Theorem \ref{new2.1},
suppose there is a unital trace preserving   completely positive map $\Phi$ satisfying {\rm (a)},
and $D$ is a doubly stochastic matrix satisfying {\rm (b)}. Let  $D$ be expressed as in (\ref{ds}).
Then the mixed unitary map $\Psi: M_n\rightarrow M_n$ defined by
$$\Psi(X) = \sum_{j=1}^r t_j (UP_jV)^*X(UP_jV)$$
also satisfies $\Psi(A_i) = B_i$ for $i = 1, \dots, k$.
\end{corollary}

\begin{remark} \label{rmk2.2} \rm The following remarks concerning Theorem \ref{new2.1} are in order.

\begin{itemize}
\item [a)] Note that the same conclusion of Theorem \ref{new2.1} holds
for completely positive maps without the unital/trace preserving requirement if the matrices $U$ and $V$
in the hypothesis are merely invertible instead of unitary.
In other words, the result applies to two families
$\{A_1, \dots, A_k\} \subseteq H_n$ and $\{B_1, \dots, B_k\}\in H_m$
such that each family is simultaneously congruent to diagonal matrices.
For example, two Hermitian matrices $A_1$ and $A_2$ are simultaneously
congruent to diagonal matrices if any one of the following conditions is satisfied.

\medskip
\hskip 1in
{\rm (i)} $\alpha_1A_1+\alpha_2A_2$ is positive definite.

\medskip
\hskip 1in
{\rm (ii)} Both $A_1$ and $A_2$ are positive semidefinite.

\medskip
\item [b)] To check conditions (b) and the corresponding ones in (1) , (2) and (3), one can use standard linear programming
techniques. In (1), one can divide the problem of finding a column stochastic $D$ into $m$ independent
problems of finding nonnegative vectors ${\mathbf d}_q = (d_{1q}, \dots, d_{nq})^t$
with entries summing up to one such that
$(a_{ij}){\bf d}_q$ equals to the $q$th column of $(b_{ij})$ for $q =1, \dots, m$.
Of course, the solution set of each of this problem is the convex
polyhedron
$$P_q = \left\{ (d_{1q}, \dots, d_{nq})^t: d_{pq} \ge 0 , \sum_{p=1}^n d_{pq}=1, \hbox{ and }
b_{iq} = \sum_{\ell=1}^n a_{i\ell}d_{\ell q}, i = 1, \dots, k
\right\}$$
in $\IR^n$. The extreme points of the polyhedron $P_q$, if non-empty,
has at most $k$ (respectively, $k+1$) nonzero entries
because we need $n$ equalities among the inequality and equality
constraints to determine an extreme point of $P_q$.
Thus, at least $n-k$ of the inequality constraints $d_{pq} \ge 0$ have to be equalities
to determine an extreme point.
By the above discussion, one can construct a sparse matrix $D$ as a solution using the
extreme points in the solution sets $P_q$.
Similarly, for (2) we can construct an extreme point (with sparse pattern) of the set
of row stochastic matrices $D$ satisfying $(b_{ij}) = (a_{ij}) D$. However, unlike (1), we cannot
treat individual rows separately to reduce the complexity of the computation.

\item [c)] In the construction of $F_1, \dots, F_r$ in the last assertion of
the theorem, we see that each $F_j$ has rank at most $\ell = \min\{m,n\}$
so that the corresponding completely positive map $\Phi$ is a super $\ell$-positive
map \cite{SSZ}.

\item [d)] Note that for the diagonal matrix $X_i = U^*A_iU$,
the map $X_i \mapsto F_j^*X_iF_j$ has only one nonzero entry
at the $(j,j)$ position obtained by taking a nonnegative
combination of the diagonal entries of $X_i$.
In the context of quantum information science, it is easy to
implement the map (quantum operation) $\Phi$, and all the actions
only take place at the diagonal entries (classical channels)
once $A_1, \dots, A_k$ and $B_1, \dots, B_k$ are in diagonal forms.
\end{itemize}
\end{remark}

\section{Completely positive maps on a single matrix}

For a single matrix, we can give a more detailed analysis of the result in Theorem \ref{new2.1},
and show that the study is related to other topics such as eigenvalue inequalities.
Moreover, the results show that
there are some results on trace preserving completely positive maps
with no analogs for unital completely positive maps, and vice versa.

\subsection{Unital completely positive maps}

\begin{theorem} \label{3.1}
Suppose  $A \in H_n$ and $B \in H_m$ have eigenvalues $a_1,\dots, a_n$ and $b_1,\dots,b_m$,
respectively. Let $\a=(a_1,\dots, a_n)$ and $\b=(b_1,\dots,b_m)$.
The following conditions are equivalent.
\begin{itemize}
\item[{\rm (a)}] There is  a (unital) completely positive map
$\Phi: M_n \rightarrow M_m$ such that $\Phi(A) = B$.
\item[{\rm (b)}] There is a nonnegative (column stochastic) matrix $D = (d_{pq})$
such that $\b = \a D$.
\item[{\rm (c)}]  There are
real numbers $\gamma_1, \gamma_2 \ge 0$ (with $\gamma_1 = \gamma_2 = 1$)
such that
\begin{equation} \label{eq-1}
\gamma_2 \min\{a_i:1\le i\le n\} \le
b_j \le \gamma_1 \max\{a_i:1\le i\le n\}.
\end{equation}
for all $1\le j\le m$.
\end{itemize}
\end{theorem}

\it Proof. \rm
By Theorem \ref{new2.1}, (a) and (b) are equivalent. Without loss of generality, we may assume that $a_1\ge\cdots\ge a_n$ and $b_1\ge\cdots\ge b_m$.

(b) $\Ra$ (c): Suppose there is a nonnegative   matrix $D = (d_{pq})$
such that $\b = \a D$.  Let $\gamma_1=\sum_{p=1}^nd_{p1}$ and $\gamma_2=\sum_{p=1}^nd_{pm}$.
Then for each $1\le j\le m$, we have
$$\gamma_2 a_n =  \sum_{p=1}^nd_{pm}a_n
\le \sum_{p=1}^nd_{pm}a_p =b_m
\le b_j \le b_1=\sum_{p=1}^nd_{p1} a_p
\le\sum_{p=1}^nd_{p1} a_1  =\gamma_1  a_1. $$
If $D$ is column stochastic, then it follows from definition that $\gamma_1 = \gamma_2 = 1$.

(c) $\Ra$ (b): Suppose (c) holds. Then for each $i=1,\dots m$, there exists $0\le t_i\le 1$ such that  $b_i=t_i\gamma_1a_1+(1-t_i)\gamma_2a_n$. Let $D$ be the  $n\times m$  matrix
$$ \[
\begin{array}{cccc}
t_1\gamma_1&t_2\gamma_1&\cdots& t_m\gamma_1\\
0&0&\cdots& 0\\
\vdots&\vdots&\vdots& \vdots\\
0&0&\cdots& 0\\
(1-t_1)\gamma_2&(1-t_2)\gamma_2&\cdots&   (1-t_m)\gamma_2\end{array}\]$$  Then $\b = \a D$ and $D$ is column stochastic if $\gamma_1=\gamma_2=1$.
\qed

\begin{corollary} \label{new2.4}
If $A\in H_n$ and $B\in H_m$ are nonzero positive semi-definite, then
there is a completely positive map $\Phi: M_n\rightarrow M_m$
such that $\Phi(A) = B$, and there is a completely positive map
$\Psi: M_m \rightarrow M_n$ such that $\Psi(B) = A$.
\end{corollary}

Using Theorem \ref{3.1}, one can construct a pair of
density matrices $(A,B) \in H_n\times H_m$
such that there is a unital completely positive map $\Phi$ such that
$\Phi(A) = B$, but there is no unital completely positive map
$\Psi$ such that $\Psi(B) = A$.

\subsection{Trace preserving completely positive maps}

\begin{theorem} \label{3.2}
Suppose  $A \in H_n$ and $B \in H_m$ have eigenvalues $a_1,\dots, a_n$
and $b_1,\dots,b_m$ respectively. Let $\a=(a_1,\dots, a_n)$ and $\b=(b_1,\dots,b_m)$.
The following conditions are equivalent.
\begin{itemize}
\item[{\rm (a)}] There is  a trace preserving completely positive map
$\Phi: M_n \rightarrow M_m$ such that $\Phi(A) = B$.
\item[{\rm (b)}] There exists an $n\times m$ row stochastic
matrix $D$  such that $\b=\a D$. Moreover, we can
assume that the $p${\rm th} and $q${\rm th} row are identical whenever
$a_pa_q > 0$, and the $p${\rm th} row of $D$ can be  arbitrary
(nonnegative with entries summing up to 1) if
$a_p = 0$.
\item[{\rm (c)}] We have $\tr A = \tr B$ and
$\sum_{p=1}^n |a_p| \ge
\sum_{q=1}^m |b_q|$.
Equivalently, $\tr A = \tr B$ and the sum of the positive
(negative) eigenvalues of $A$ is not smaller (not larger)
than the sum of  the positive (negative) eigenvalues of $B$.
\end{itemize}
\end{theorem}

\it Proof. \rm For simplicity, we assume that $a_1 \ge \cdots a_r\ge 0 > a_{r+1}
\ge \cdots \ge a_n$ and $b_1\ge\cdots\ge b_s \ge 0> b_{s+1}\ge\cdots\ge b_m$. Let $a_+=\sum_{p=1}^ra_p$, $a_-=\sum_{p=r+1}^na_p$ and $b_+=\sum_{q=1}^sb_q$, $b_-=\sum_{q=s+1}^mb_q$.

By Theorem \ref{new2.1}, we have (b) $\Ra$ (a).
Also, by Theorem \ref{new2.1},  if (a) holds, then $\b = \a D$ for an $n\times m$ row
stochastic matrix $D$.
Next, we show that $D$ can be chosen to satisfy the second assertion of condition (b).
To this end, let $\d_1, \dots, \d_n$ be the rows of $D$. Clearly, if $a_p = 0$, we can replace
the $p\,$th row of $D$ by any nonnegative vectors with entries summing up to 1 to get
$\tilde D$ and we still have $\a \tilde D = \a D = \b$.
Now, suppose $a_1 > 0$. Then $a_+ > 0$.
We can replace the first $r$ rows (or the rows correspond to $a_p > 0$) by
$$\d_+ = \sum_{j=1}^r \frac{a_j}{a_+} \d_j$$
to obtain $\tilde D$.
Then $\d_+$ has nonnegative entries summing up to 1, and
$\a \tilde D = \a D = \b$. Similarly, suppose  $a_n < 0$.
Then $a_- < 0$. We can further replace the last $n-s$ row of $D$ by
$$\d_- = \sum_{j=s+1}^n \frac{a_j}{a_-} \d_j$$
to obtain $\tilde D$.
Then $\d_-$ has nonzero entries summing up to 1, and
$\a \tilde D = \a D = \b$.

(b) $\Ra$ (c): Suppose $\b=\a D$. We have
$$\tr B=\sum_{q=1}^mb_q=\sum_{q=1}^m\sum_{p=1}^na_pd_{pq}
=\sum_{p=1}^na_p\(\sum_{q=1}^md_{pq}\)=\sum_{p=1}^na_p=\tr A$$
and
$$
\sum_{q=1}^m|b_q|=\sum_{q=1}^m| \sum_{p=1}^na_pd_{pq}|
\le \sum_{q=1}^m \sum_{p=1}^n|a_p|d_{pq}=\sum_{p=1}^n|a_p|\(\sum_{q=1}^md_{pq}\)=\sum_{p=1}^n|a_p|\,.
$$
Since
$$\sum_{i=1}^n a_i = \sum_{j=1}^m b_j, \quad
\sum_{i=1}^n |a_i| = a_+ - a_- = \sum_{i=1}^n a_i - 2 a_- = 2 a_+ - \sum_{i=1}^n a_i,$$
and
$$\sum_{j=1}^m |b_j| = b_+ - b_- = \sum_{j=1}^m b_j - 2 b_- = 2 b_+ - \sum_{j=1}^m b_j,$$
the last assertion of (c) follows.

(c) $\Ra$ (b):
Suppose $\tr A = \tr B$, $a_+ \ge b_+$ and $a_-\le b_-$.
\iffalse
Since $a_++a_-=\tr A=\tr B=b_++b_-$ and $a_+-a_-=\sum_{p=1}^n |a_p| \ge
\sum_{q=1}^m |b_q|=b_+-b_-$, we have $a_+\ge b_+$ and $a_-\le b_-$.
\fi
Let $$t_q=\left\{\begin{array}{ll} \dfrac{b_q}{a_+} &\mbox{ for } 1\le q\le s, \\&\\
\dfrac{b_q}{a_-}&\mbox{ for } s < q\le m .\end{array}\right.$$
Here, if $a_+ = 0$ then $b_+ = 0$, and we can set $t_q = 0$ for $1 \le q \le s$.
If $a_- = 0$ then $b_- = 0$, and  $s = m$.
Therefore, $t_q\ge 0$ for all $1\le q\le m$. We have
$$
a_+ \ge b_+ = (a_+) \sum_{q=1}^{s}t_q,
\quad \mbox{ and } \quad |a_-| \ge |b_-| = |a_-|\sum_{q=s+1}^{n}t_q.$$
Let $u=1-\sum_{q=1}^{s}t_q\ge 0$, $v=1-\sum_{q=s+1}^{n}t_q\ge 0$ and
$ D $ be an $n\times m$ row stochastic matrix with
$$
  p\,\mbox{th row }=\left\{\begin{array}{ll}  (t_1,t_2,\dots,t_s,0,\dots,0,u) &\mbox{ for }1\le p\le r,
  \\&\\
  (0,\dots,0,t_{s+1},t_{s+2},\dots,t_{m-1},t_m+v)&\mbox{ for }r+1\le p\le n\,.\end{array}\right.$$
Since $$ua_++va_-=(a_+-b_+)+(a_--b_-)=0,$$
we have $\b=\a D$. \qed

\begin{corollary} \label{new3.4}
If $A\in H_n$ and $B\in H_m$ are density matrices, i.e. positive semi-definite and $\tr A=\tr B=1$, then
there is a completely positive map $\Phi: M_n\rightarrow M_m$
such that $\Phi(A) = B$, and there is a completely positive map
$\Psi: M_m \rightarrow M_n$ such that $\Psi(B) = A$.
\end{corollary}

As remarked after Corolary \ref{new2.4},
one may not be able to find a unital completely positive
map taking a density matrix  $B \in H_m$ to another
density matrix $A \in H_n$ even if there is a unital
positive completely positive map sending $A$ to $B$.

\subsection{Unital trace preserving completely positive maps}

 \

Suppose there is a unital completely positive map
sending $A$ to $B$, and also a trace preserving completely positive
map sending $A$ to $B$. Is there a unital trace preserving
completely positive map sending $A$ to $B$?
The following example shows that the answer is negative.

\begin{example} \label{2.6}
Suppose $A = \diag(4,1,1,0)$ and $B = \diag(3,3,0,0)$.
By Theorems \ref{3.1} and \ref{3.2} there is a trace preserving
completely positive map sending $A$ to $B$, and also a unital
completely positive map sending $A$ to $B$. Let
$A_1=A-I_4=\diag(3,0,0,-1)$ and
$B_1=B-I_4=\diag(2,2,-1,-1)$. By Theorem \ref{3.2},  there is no
trace preserving  completely positive map sending
$A_1$ to $B_1$. Hence, there is no  unital trace preserving
completely positive map sending $A$ to $B$.
\end{example}

As shown in Corollary \ref{new2.2},
if there is a unital trace preserving map sending a commuting family in $H_n$
to a commuting family in $H_m$, then we may chose the map to be mixed unitary.
In the following, we show that for the case when $k = 1$ in Theorem \ref{new2.1}, one can even assume that the map is the average
of $n$ unitary similarity transforms.

Let $\a=(a_1, \dots, a_n), \b=(b_1, \dots, b_n)\in\IR^n$. We say that $\b$ is majorized by $\a$ ($\b\prec \a$) if for every $1\le k< n$, the   sum of the $k$ largest entries of $\b$ is less than or equal to  the sum of the $k$ largest entries of $\a$, and $\sum_{i=1}^na_i=\sum_{i=1}^nb_i$.

\begin{theorem} \label{3.5}Suppose  $A \in H_n$ and $B \in H_m$ have eigenvalues $a_1,\dots, a_n$ and $b_1,\dots,b_m$ respectively. Let $\a=(a_1,\dots, a_n)$ and $\b=(b_1,\dots,b_m)$.
The following are equivalent.
\begin{enumerate}

\item[{\rm (a)}] There exists a unital trace preserving completely positive
map $\Phi$ such that $\Phi(A)=B$.

\item[{\rm (a1)}] There exists a mixed unitary  completely positive
map $\Phi$ such that $\Phi(A)=B$.

\item[{\rm (a2)}]  There exist unitary matrices  $U_1, \dots, U_n \in M_n$ such that
$B = \frac{1}{n}\sum_{j=1}^n U_j^* A U_j$.

\item[{\rm (a3)}] For each $t\in \IR$, there exists
a  trace preserving completely positive
map $\Phi_t$ such that $\Phi_t(A-tI)=B-tI$.

\item[{\rm (b)}] There is a doubly stochastic matrix $D$ such that
$\b = \a D$.

\item[{\rm (b1)}] There is a unitary matrix $W\in M_n$ such that
$W\diag(a_1, \dots, a_n)W^*$ has diagonal entries $b_1, \dots, b_n$.

\item[{\rm (c)}] $\b \prec \a$.

\end{enumerate}
Moreover, if condition (b) holds and $D = \sum_{\ell=1}^r t_\ell P_\ell$
such that $t_1, \dots, t_r$ are positive numbers summing up to 1
and  $P_1, \dots, P_r$ are permutation matrices, then
$B = \sum_{j=1}^r V^*P_j^tU^*AUP_j V$, where
$U^*AU = \diag(a_1, \dots, a_n)$ and $V^*BV = \diag(b_1, \dots, b_n)$.
\end{theorem}

\it Proof. \rm (b) $\iff$ (c) is a standard result of majorization; see \cite{MO}.

The implications
(a2) $\Ra $ (a1) $\Ra$ (a) $\Ra$ (a3) are obvious.

(a3) $\Ra$ (c) :  We may assume that $a_1 \ge \cdots \ge a_n$
and $b_1 \ge \cdots \ge b_n$.
For $1\le k<n$ choose $t$ such that  $a_k\ge t\ge a_{k+1}$.
Then there is a trace preserving completely positive linear map
$\Phi_t$ such that $\Phi_t(A-tI_n) = B-tI_n$.
By Theorem \ref{3.2}, the sum of the $k$ positive eigenvalues of $B-tI_n$
is no larger than that of $A-tI_n$. Thus,
\begin{equation}\label{eq1}
\sum_{i=1}^k b_i - kt = \sum_{i=1}^k (b_i - t)
\le \sum_{i=1}^k (a_i - t) = \sum_{i=1}^k a_i - kt \, .
\end{equation}
We see that $\sum_{i=1}^k b_i \le \sum_{i=1}^k a_i$ for
$k = 1, \dots, n-1$. Since $\tr A = \tr B$, we have $\b \prec \a$.

(c) $\Ra$ (b1) is a result of Horn \cite{H}.

(b1) $\Ra$ (a2) :   Suppose $W$ is a unitary matrix such that
the diagonal of $W\diag(a_1, \dots, a_n)W^*$ has diagonal entries
$b_1, \dots, b_n$. Let $w = e^{i2\pi/n}$,
$P = \diag(1, w, \dots, w^{n-1})$. Then
$$B = \frac{1}{n} \sum_{j=1}^n UW(P^j)^*WU A U^*W^*(P^j)U^*.$$
Thus, (a2) holds.
\qed

In the context of quantum information theory, a completely positive
map in condition (a1) of the above theorem is a {\it mixed unitary
quantum channel/operation}. By the above theorem,  the existence of a mixed
unitary quantum channel $\Phi$ taking a quantum state $A\in H_n$
to a quantum state $B \in H_m$ can be described in terms of trace
preserving completely positive maps, namely, condition (a3).
However, despite the duality of
the two classes of maps, there is no  analogous condition
in terms of unital completely positive map; see Proposition \ref{1.2}.

\section{Additional remarks and future research}

To study unital completely positive maps connecting two families
$\{A_1, \dots, A_k\} \subseteq H_n$ and $\{B_1, \dots, B_k\} \subseteq H_m$,
one can use the results on completely
positive maps and add $I_n$ and $I_m$ to the two families.

The following result shows that the study of a completely positive
maps sending $A_1, \dots, A_k \in H_n$ to $B_1, \dots, B_k \in H_m$
can be reduced to the study of unital completely positive maps.

\begin{theorem} Let $A_1, \dots, A_k \in H_n$, and $B_1, \dots, B_k \in H_m$.
There is a completely positive map $\Phi: M_n \rightarrow M_m$ such that
$\Phi(A_j) = B_j$ for $j = 1, \dots, k$  if and only if there exists $\gamma >0$ and  a unital
completely positive map $\Psi: M_{n+1} \rightarrow M_m$ such that
$\Psi(A_j \oplus [0]) = \gamma^{-1} B_j$ for $j = 1, \dots, k$.
\end{theorem}

\it Proof. \rm
Suppose $\Phi: M_n \rightarrow M_m$ has operator sum representation
(\ref{operatorsum}) and
satisfies $\Phi(A_j) = B_j$ for $j = 1, \dots, k$.
Let $\Phi(I_n) = P \in H_m$.  Choose $\gamma >0$ such that   $\gamma I_m - P$ is positive semi-definite. Then we have $\gamma I_m - P=\sum_{j=1}^s g_j ^*g_j$ for some $1\times m$ matrices $g_j $.
For $j =1, \dots, \tilde r$ with $\tilde r = \max\{r,s\}$,
let $\tilde F_j\in M_{n+1,m}$ be such that
$\tilde F_j = \[\begin{array}{c}F_j\\  g_j\end{array}\]$, where $F_j = 0$ if $j > r$
and $g_j = 0$ if $j > s$.
Define $\Psi: M_{n+1} \rightarrow M_m$ by
$\Psi(X) =\dfrac{1}{\gamma} \sum_{j=1}^{\tilde r} \tilde F_j^* X \tilde F_j$.
One readily checks that $\Psi(I_{n+1}) = I_m$
and $\Psi(A_j\oplus [0]) = \gamma^{-1} B_j$ for $j = 1, \dots, k$.

Conversely, suppose $\gamma>0$ and  $\Psi: M_{n+1} \rightarrow M_m$ is a unital
completely positive
map such that $\Psi(A_j\oplus [0]) = \gamma^{-1} B_j$ for $j =1, \dots, k$,
then one can check $\Phi:M_n \rightarrow M_m$
defined by
$\Phi(X) = \gamma \Psi(X\oplus [0])$ is a completely positive map
satisfying $\Phi(A_j) = B_j$ for $j = 1, \dots, k$.
\qed

\medskip

Finding a unital completely positive map connecting two
general (non-commuting) families of Hermitian matrices is very challenging.
In the case of non-commuting families $\{A_1, A_2\} \subset H_n$,
and $\{B_1,B_2\} \subseteq H_m$ with two elements,
the problems reduce to the study of unital completely positive maps $\Phi$
satisfying $\Phi(A_1+iA_2)=B_1+iB_2$. For $n = 2, 3$. There are
partial answers of the problem
in terms of the {\it numerical range} and {\it dilation} of operators.
Recall that the numerical range of $T \in M_n$ is the set
$$W(T) = \{ x^*Tx: x \in \IC^n, \ x^*x = 1\},$$
and $T$ has a {\it dilation} $\tilde T \in M_m$ if there is an $n\times m$
matrix $X$ such that $XX^* = I_n$ and $X\tilde T X^* = T$.
We have the following result; see \cite{CL1,CL2}.

\begin{theorem} \label{4.2} Let $(A,B) \in M_n\times M_m$. Suppose $n = 2$,
or $n=3$ such that $A$ is unitarily reducible, i.e., $A$ is untiarily
similar to $A_1 \oplus [\alpha]$ for some $A_1 \in M_2$ and $\alpha \in \IC$.
Then the following conditions are equivalent.
\begin{itemize}
\item[{\rm (a)}]
There is a unital completely positive map $\Phi: M_n \rightarrow M_m$
such that $\Phi(A) = B$.
\item[{\rm (b)}]
$B$ has a dilation of the form $I_r \otimes A$.
\item[{\rm (c)}]
$W(B) \subseteq W(A)$.
\end{itemize}
\end{theorem}

Special cases of the above theorem include the case when $A\in M_3$ is a normal matrix.
However,
there are examples showing that the result fails if $A$ is an arbitrary matrix in $M_3$
or an arbitrary normal matrix in $M_4$; see \cite{CL1}.

\medskip
In connection to Theorem \ref{4.2}, one may ask for the condition of $A\in M_n$ to be a
dilation of $B\in M_m$ itself. The problem is challenging even for normal
matrices $A$ and $B$; see \cite{Ho}.

\medskip
Also, it is interesting to impose condition on the Kraus (Choi) rank, i.e., the
minimum number of matrices $F_1, \dots, F_r$ needed in the operator sum representation
of the completely positive maps. As mentioned in Proposition \ref{1.3},
the study is related to the study of principal submatrices of a Hermitian matrices,
which is related to the study of spectral inequalities and Littlewood-Richardson
rule; see \cite{FFLP,LP03}.

\bigskip\noindent
{\large{\bf Acknowledgment}}

Part of the results in the paper was reported at the Workshop on
Mathematics in Experimental Quantum Information Processing, IQC,
Waterloo, August 10-14, 2009. The comments of the participants
and the support of the organizer are graciously acknowledged.
Research of both authors are supported by USA NSF. The first author
was also supported by a HK RCG grant, and the
Key Disciplines of Shanghai Municipality Grant S30104.
He is an honorary professor of the University of Hong Kong,
the Taiyuan University of Technology,
and the Shanghai University.

\bibliographystyle{amsplain}


\end{document}